\documentclass[12pt]{article}
\usepackage{amssymb,bm}
\usepackage{amsmath}
\usepackage{amsthm}
\usepackage{graphicx}
\usepackage{hyperref} 
\hypersetup{pdfborder=0 0 0}
\newtheorem{theorem}{{\sc Theorem}}[section]

\bmdefine\Bone{1}

\def\XXint#1#2#3{{\setbox0=\hbox{$#1{#2#3}{\int}$ }
\vcenter{\hbox{$#2#3$ }}\kern-.6\wd0}}

\newcommand{\re}{\Re\mathfrak{e}}

\newcommand{\Ge}{\epsilon}

\bmdefine\BGa{\alpha}
\bmdefine\BGb{\beta}
\bmdefine\BGd{\delta}
\bmdefine\BGe{\epsilon}
\bmdefine\BGve{\varepsilon}
\bmdefine\BGf{\phi}
\bmdefine\BGvf{\varphi}
\bmdefine\BGg{\gamma}
\bmdefine\BGc{\chi}
\bmdefine\BGi{\iota}
\bmdefine\BGk{\kappa}
\bmdefine\BGl{\lambda}
\bmdefine\BGn{\eta}
\bmdefine\BGm{\mu}
\bmdefine\BGv{\nu}
\bmdefine\BGp{\pi}
\bmdefine\BGth{\theta}
\bmdefine\BGvth{\vartheta}
\bmdefine\BGr{\rho}
\bmdefine\BGvr{\varrho}
\bmdefine\BGs{\sigma}
\bmdefine\BGvs{\varsigma}
\bmdefine\BGt{\tau}
\bmdefine\BGj{\tau}
\bmdefine\BGu{\upsilon}
\bmdefine\BGo{\omega}
\bmdefine\BGx{\xi}
\bmdefine\BGy{\psi}
\bmdefine\BGz{\zeta}
\bmdefine\BGD{\Delta}
\bmdefine\BGF{\Phi}
\bmdefine\BGG{\Gamma}
\bmdefine\BGL{\Lambda}
\bmdefine\BGP{\Pi}
\bmdefine\BGT{\Theta}
\bmdefine\BGS{\Sigma}
\bmdefine\BGU{\Upsilon}
\bmdefine\BGO{\Omega}
\bmdefine\BGX{\Xi}
\bmdefine\BGY{\Psi}

\bmdefine\BFM{\mathfrak{M}}
\bmdefine\BFb{\mathfrak{b}}
\bmdefine\BFk{\mathfrak{k}}
\bmdefine\BFm{\mathfrak{m}}
\bmdefine\BFu{\mathfrak{u}}
\bmdefine\BFv{\mathfrak{v}}

\newcommand{\CH}{{\mathcal H}}

\newcommand{\CK}{{\mathcal K}}

\bmdefine\BCA{{\mathcal A}}
\bmdefine\BCB{{\mathcal B}}
\bmdefine\BCC{{\mathcal C}}
\bmdefine\BCD{{\mathcal D}}
\bmdefine\BCE{{\mathcal E}}
\bmdefine\BCF{{\mathcal F}}
\bmdefine\BCG{{\mathcal G}}
\bmdefine\BCH{{\mathcal H}}
\bmdefine\BCI{{\mathcal I}}
\bmdefine\BCJ{{\mathcal J}}
\bmdefine\BCK{{\mathcal K}}
\bmdefine\BCL{{\mathcal L}}
\bmdefine\BCM{{\mathcal M}}
\bmdefine\BCN{{\mathcal N}}
\bmdefine\BCO{{\mathcal O}}
\bmdefine\BCP{{\mathcal P}}
\bmdefine\BCQ{{\mathcal Q}}
\bmdefine\BCR{{\mathcal R}}
\bmdefine\BCS{{\mathcal S}}
\bmdefine\BCT{{\mathcal T}}
\bmdefine\BCU{{\mathcal U}}
\bmdefine\BCV{{\mathcal V}}
\bmdefine\BCW{{\mathcal W}}
\bmdefine\BCX{{\mathcal X}}
\bmdefine\BCY{{\mathcal Y}}
\bmdefine\BCZ{{\mathcal Z}}

\bmdefine\Bzr{ 0}
\bmdefine\Ba{ a}
\bmdefine\Bb{ b}
\bmdefine\Bc{ c}
\bmdefine\Bd{ d}
\bmdefine\Be{ e}
\bmdefine\Bf{ f}
\bmdefine\Bg{ g}
\bmdefine\Bh{ h}
\bmdefine\Bi{ i}
\bmdefine\Bj{ j}
\bmdefine\Bk{ k}
\bmdefine\Bl{ l}
\bmdefine\Bm{ m}
\bmdefine\Bn{ n}
\bmdefine\Bo{ o}
\bmdefine\Bp{ p}
\bmdefine\Bq{ q}
\bmdefine\Br{ r}
\bmdefine\Bs{ s}
\bmdefine\Bt{ t}
\bmdefine\Bu{ u}
\bmdefine\Bv{ v}
\bmdefine\Bw{ w}
\bmdefine\Bx{ x}
\bmdefine\By{ y}
\bmdefine\Bz{ z}
\bmdefine\BA{ A}
\bmdefine\BB{ B}
\bmdefine\BC{ C}
\bmdefine\BD{ D}
\bmdefine\BE{ E}
\bmdefine\BF{ F}
\bmdefine\BG{ G}
\bmdefine\BH{ H}
\bmdefine\BI{ I}
\bmdefine\BJ{ J}
\bmdefine\BK{ K}
\bmdefine\BL{ L}
\bmdefine\BM{ M}
\bmdefine\BN{ N}
\bmdefine\BO{ O}
\bmdefine\BP{ P}
\bmdefine\BQ{ Q}
\bmdefine\BR{ R}
\bmdefine\BS{ S}
\bmdefine\BT{ T}
\bmdefine\BU{ U}
\bmdefine\BV{ V}
\bmdefine\BW{ W}
\bmdefine\BX{ X}
\bmdefine\BY{ Y}
\bmdefine\BZ{ Z}

\DeclareMathAlphabet{\pazocal}{OMS}{zplm}{m}{n}
\newcommand{\PH}{\pazocal{H}}

\newcommand{\PS}{\pazocal{S}}
\newcommand{\PK}{\pazocal{K}}
\newcommand{\PR}{\pazocal{R}}

\DeclareMathAlphabet{\cg}{OMS}{zplm}{m}{n}

\usepackage{xcolor}

\begin{document}

\title{On the optimal analytic continuation from discrete data}
  
\author{Narek Hovsepyan}
\date{}
\maketitle

\begin{abstract}
We consider analytic functions from a reproducing kernel Hilbert space. Given that such a function is of order $\epsilon$ on a set of discrete data points, relative to its global size, we ask how large can it be at a fixed point outside of the data set. We obtain optimal bounds on this error of analytic continuation and describe its asymptotic behavior in $\epsilon$. We also describe the maximizer function attaining the optimal error in terms of the resolvent of a positive semidefinite, self-adjoint and finite rank operator.         
\end{abstract}

\section{Introduction}

Analytic functions are of central importance in many applied problems. They appear in diverse areas, such as medical imaging \cite{epst08}, nuclear physics \cite{capr74,capr79} and optimal design problems \cite{lipt01a,lipt01}. For example Fourier (or Laplace) transforms of real-valued functions vanishing on negative semi-axis are analytic in the upper half-plane. Such functions describe linear, time-invariant and causal systems. Some concrete examples are the complex magnetic permeability and complex dielectric permittivity functions \cite{lali60:8,feyn64}, complex impedance and admittance functions of electrical circuits \cite{brune31}. Further examples include transfer functions of digital filters \cite{grs01}, the dependence of effective moduli of composites on the moduli of its constituents \cite{berg78,mi80} etc.

Typically these functions can be measured only on a subset $\Gamma$ of their domain of analyticity (or its boundary). During the measurement process unavoidable error occurs, as a result an analytic function is known on $\Gamma$ up to a certain precision of order $\epsilon>0$. In order to predict the behavior of the system and to expand its application horizon, often times one is interested in extrapolating from the measured data to a given point $z$ lying outside of the data set $\Gamma$. On one hand, working in the class of analytic functions we expect rigidity (in the sense that values of an analytic function on $\Gamma$ affect its values elsewhere). On the other hand, recent work \cite{deto18, trefe19, grho-annulus, grho-gen} shows, that surprisingly there is also flexibility, meaning that the measured data can be matched up to the given precision $\epsilon$ by two analytic functions that are very different outside of the data set. 

Let $\PS=\PS(\Omega)$ denote a class of physically admissible functions that are analytic in a domain $\Omega$ of the complex plane. Aside from analyticity, the set $\PS$ may also contain further physical restrictions (cf. \cite{grho-CEMP}), such as certain symmetry constraints, asymptotic constraints at infinity or inequality constraints e.g. nonnegative imaginary parts (which can be interpreted as presence of energy loss in the system). Let $\Gamma \subset \Omega$ denote the data set where the measurements are done with relative error $\Ge$, with respect to some norm $\|\cdot\|_\Gamma$ on $\Gamma$. To quantify the flexibility of the class $S$ we ask the following question: given two functions from $\PS$ that are $\epsilon$-close on $\Gamma$ (relative to their total size on $\Omega$), how much can they differ at a point $z \in \Omega \backslash \Gamma$? Assume that the total size of a function on $\Omega$ is measured in some norm $\|\cdot\|_\Omega$, then we arrive at the quantity

\begin{equation} \label{Delta abstract}
\Delta_z(\epsilon) = \sup \left\{ \frac{|\phi(z)-\psi(z)|}{\max(\|\phi\|_\Omega, \|\psi\|_\Omega)}: \phi,\psi \in \PS \quad \text{with} \quad \frac{\|\phi-\psi\|_\Gamma}{\max(\|\phi\|_\Omega, \|\psi\|_\Omega)} \leq \epsilon \right\}.
\end{equation}

Another related quantity interesting in its own right, is the relative error of analytic continuation. Loosely speaking \footnote{For a rigorous comparison of quantities \eqref{Delta abstract} and \eqref{A abstract} in the context of the complex electromagnetic permittivity function we refer to \cite{grho-CEMP} (in this case $\PS$ is a cone of functions that is related to Herglotz-Nevanlinna functions)}, consider the difference $f = \phi - \psi$ and rescale it (say $\PS$ is a cone) by the total norm of $\phi$ and $\psi$ on $\Omega$, assume also that $\PS-\PS$ can be approximated by functions from some normed space $\PH = \PH(\Omega)$ of analytic functions in $\Omega$. We arrive at an analogous question: given that $f \in \PH$ is of order $\epsilon$ on $\Gamma$ and is of order $1$ on $\Omega$, how large can it be at the point $z$? So to quantify the stability of analytic continuation in the normed space $\PH$ we introduce

\begin{equation} \label{A abstract}
A_z(\epsilon) = \sup \left\{ |f(z)|: f \in \PH \quad \text{with} \quad \|f\|_\PH \leq 1, \ \|f\|_\Gamma \leq \epsilon  \right\}.
\end{equation}

In the setting of Hilbert spaces and when $\Gamma \Subset \Omega$ is a curve with $\|\cdot\|_\Gamma$ denoting the $L^2(\Gamma)$-norm (with respect to the arclength measure) we analyzed \eqref{A abstract} in \cite{grho-annulus, grho-gen}, where we derived optimal bounds for it and showed that it behaves like a power law: $A_z(\epsilon) \approx \epsilon^{\gamma(z)}$, where the exponent $\gamma(z) \in (0,1)$ decreases to $0$, as we move further away from the source of
data. How fast $\gamma(z)$ decays depends strongly on the geometry of the domain and the data source. The most common setting, where \eqref{A abstract} is analyzed in the literature is in the space
of bounded analytic functions. The power law estimates are then derived from a maximum modulus principle, the Hadamard three-circles theorem is a classical example of such estimate. For related works we refer the reader to \cite{davis52,ciulli69,mill70,payne75,fran90,vese99,fdfd09,deto18,
trefe19}.

This paper is dedicated to the analysis of \eqref{A abstract} in the Hilbert space setting, when $\Gamma = \{z_j\}_{1}^n$ represents a finite set of distinct points, where the function values are measured. In this case $\|\cdot\|_\Gamma$ is a seminorm, so we use the notation $[\cdot]_\Gamma$ instead, and treating all the points equally we consider the $l^2$-seminorm: $[f]^2_\Gamma = \sum_j |f(z_j)|^2$. The first difference of the discrete setting vs. the continuum one is that in the former case an analytic function is not determined uniquely by its values on $\Gamma$, as a result $A_z(\epsilon)$ does not converge to zero as $\epsilon \to 0$. So then the questions are what is $A_z(0)$ and what is the next term in the asymptotic expansion of $A_z(\epsilon)$. The answer to the last question reveals the second key distinction of the discrete setting, showing that there is no fractional power of $\epsilon$ and the correction term is of order $\epsilon$. Namely, we will characterize $A_z(0)$ (in terms of the reproducing kernel of the space $\PH$, cf. Theorem~\ref{THM discrete opt rec}) and show that

\begin{equation} \label{A_z asymp}
A_z(\epsilon) = (1 + \sigma \epsilon) A_z(0) + O(\epsilon^2),
\end{equation}

\noindent where $\sigma = \sigma(z)>0$ will depend on the space $\PH$ and the data set $\Gamma$. Note that the set of values $V_z(\epsilon) = \left\{ f(z): \|f\|_\PH \leq 1, \ [f]_\Gamma \leq \epsilon  \right\}$ is a convex, centrally symmetric ($c \in V_z$ iff $-c \in V_z$) subset of the complex plane, and $A_z(\epsilon)$ is its "radius", i.e. half of the diameter. The formula \eqref{A_z asymp} then shows the relation between the radii $A_z(0), \ A_z(\epsilon)$ of the original and perturbed function value sets, respectively. 

The quantity \eqref{A abstract} is also related to the optimal estimation of the point evaluation functional $f \mapsto f(z)$ (see \cite{micriv77, micriv85} and references therein for the general theory of optimal estimations and optimal recovery). Following \cite{micriv77} let us formulate the question of optimal recovery. Let $f_j' = f(z_j) + \delta_j$ represent the erroneous measurement of the function value at the point $z_j$ for $j=1,...,n$. Assume, that the error is of order $\epsilon$, namely let $\bm{\delta} = (\delta_1,...,\delta_n)$ be the error vector and let $|\bm{\delta}| \leq \epsilon$, where $|\cdot|$ is the Euclidean length of a vector in $\mathbb{C}^n$. The task is to approximate $f(z)$ at a fixed point $z \notin \Gamma=\{z_j\}_1^n$. The error of a linear estimation algorithm (it is enough to restrict consideration only to linear algorithms \cite{maos75}) is then defined as

\begin{equation} \label{E_z(eps,c)}
E_z(\epsilon,\bm{c}) = \sup \left\{ |f(z) - \sum_{j=1}^n c_j (f(z_j) + \delta_j)| : \ \|f\|_\PH \leq 1, \ |\bm{\delta}| \leq \epsilon  \right\},
\end{equation}  

\noindent where $\bm{c} = (c_1,...,c_n) \in \mathbb{C}^n$ is a given vector defining the linear algorithm and the supremum goes over all $f \in \PH$ and $\bm{\delta} \in \mathbb{C}^n$ satisfying the above-mentioned constraints. The intrinsic error of the estimation problem is

\begin{equation} \label{E_z(eps)}
E_z(\epsilon) = \inf_{\bm{c}} E_z(\epsilon, \bm{c}).
\end{equation}

\noindent Any algorithm achieving this infimum yields an optimal procedure for estimating $f(z)$. Theorem 1 of \cite{maos75} implies that

\begin{equation} \label{A_z=E_z}
A_z(\epsilon) = E_z(\epsilon).
\end{equation} 

\noindent Let us actually prove this equality using an idea from \cite{fisher83} (Section 7.5). The constraints in \eqref{E_z(eps,c)} are invariant under multiplying $f$ and $\bm{\delta}$ with a constant phase factor, so instead of maximizing the absolute value in \eqref{E_z(eps,c)} we can equivalently maximize the real part. Next, applying von Neumann's minimax theorem \cite{nik54} we obtain

\begin{equation*}
E_z(\epsilon, \bm{c}) = \inf_{\bm{c}} \sup_{f, \bm{\delta}} \re \left\{ f(z) - \sum_{j=1}^n c_j (f(z_j) + \delta_j) \right\} = \sup_{f, \bm{\delta}} \inf_{\bm{c}} \re \left\{ f(z) - \sum_{j=1}^n c_j (f(z_j) + \delta_j) \right\}. 
\end{equation*}

\noindent It remains to note that the inner infimum will be $-\infty$ unless $f(z_j) + \delta_j = 0$ for all $j$. This implies that the supremum can be restricted to considering those $f \in \PH$ with $\|f\|_{\PH} \leq 1$ for which the choice $\delta_j = - f(z_j)$ can be made, which means that $f$ must also satisfy the second constraint $|\bm{\delta}| = [f]_\Gamma \leq \epsilon$. This concludes the proof of \eqref{A_z=E_z}.

In \cite{mafed01} the authors analyze a quantity related to $E_z(\epsilon)$, namely in order to obtain constructive results in \eqref{E_z(eps,c)} they replace the target functional with the square root of $|f(z) - \sum_{j} c_j f(z_j)|^2 + |\sum_j c_j \delta_j)|^2$. The square of the replaced quantity is comparable to $E_z(\epsilon)$ and hence also to $A_z(\epsilon)$. In this work we take an alternative approach and analyze $A_z(\epsilon)$ directly using variational methods and derive the asymptotic expansion result \eqref{A_z asymp} that is analogous to that of \cite{mafed01} (see Theorem 4 therein). Further, we do not assume linear independence of the point evaluation functionals $f \mapsto f(z_j)$ for $j=1,...,n$. Moreover, we describe the maximizer function attaining the supremum in \eqref{A abstract} via the resolvent of a positive semidefinite, self-adjoint and finite rank operator, which (by taking limits as $\epsilon \to 0$) also allows us to obtain a characterization for $A_z(0)$.

\section{The Main Result}
\setcounter{equation}{0}

\noindent Let $\PH = \PH(\Omega)$ be a Hilbert space of analytic functions in a domain $\Omega \subset \mathbb{C}$. Assume that the point evaluation functional $f \mapsto f(z)$ is continuous for any point $z \in \Omega$, then by the Riesz representation theorem, there exists an element $p_z \in \PH$ such that $f(z) = (f,p_z)_{\PH}$. So inner products with the function $p(\zeta, z):=p_z(\zeta)$ reproduce values of a function in $\PH$. In this case $\PH$ is called a a reproducing kernel Hilbert space (RKHS) with kernel $p$ (cf. \cite{para16}). Examples of such spaces include the Hardy spaces $H^2$, the Bergman spaces $A^2$ etc. From now on let us drop the subscript $\PH$ from the notation of the norm and the inner product of $\PH$.

Let $z, z_1,...,z_n \in \Omega$ be distinct points and set $[f]^2 = \sum_j |f(z_j)|^2$. Consider the problem

\begin{equation} \label{A abstract discrete}
A_z(\epsilon) = \sup \left\{ |f(z)|: f \in \PH \quad \text{with} \quad \|f\| \leq 1, \ [f] \leq \epsilon  \right\}.
\end{equation} 

\noindent Introduce the restriction operator $\PR: \PH \to \mathbb{C}^n$ given by $\PR f = \bm{f}:= (f(z_1),...,f(z_n))$ and let $\PK = \PR^* \PR$ (cf. \cite{davis52}), then it is easy to see that

\begin{equation} \label{K}
\CK f = \sum_{j=1}^n f(z_j) p_{z_j}
\end{equation}

\noindent and $(\CK f, g) = (\bm{f}, \bm{g})_{\mathbb{C}^n}$, in particular the second constraint of \eqref{A abstract discrete} can be rewritten via the quadratic form of $\PK$:

\begin{equation} \label{K quadratic form}
(\CK f, f) = [f]^2.
\end{equation}

\noindent Clearly $\CK$ is a self-adjoint, positive semidefinite and compact (in fact, finite rank) operator on $\PH$. Set $U_0=\ker \CK$ \footnote{when $\{p_{z_j}\}$ are linearly independent $U_0 = \{f \in \CH: f(z_1)=...=f(z_n)=0\}$}, let $\{\lambda_1,...,\lambda_m\}$ be all distinct nonzero eigenvalues of $\CK$ and set $U_j = \ker (\CK - \lambda_j I), \quad j=1,...,m$. Hence, $\CH = \bigoplus_{j=0}^m U_j$, and let $P_j : \CH \to U_j$ denote the orthogonal projection onto the closed subspace $U_j$ for $j=0,...,m$.

\begin{theorem} \label{THM discrete opt rec}
With the notation introduced above and assuming $p_z \notin \PK(\PH)$, we have

\begin{equation} \label{main identity}
A_z(\epsilon) = (1 + \sigma \epsilon) A_z(0) + O(\epsilon^2),
\end{equation}

\noindent where, in fact $A_z(0) = \|P_0 p_z \|$ and

\begin{equation} \label{sigma}
\sigma^2 = \frac{1}{\|P_0 p_z \|^2} \sum_{j=1}^m \frac{\|P_j p_z\|^2}{\lambda_j}.
\end{equation}

\noindent Moreover, the maximizer function for $A_z(\epsilon)$ is given via the resolvent operator  

\begin{equation} \label{A_z using u_eps}
A_z(\epsilon) = \frac{u_\epsilon(z)}{\|u_\epsilon\|}, \qquad \qquad u_\epsilon = (\CK + \eta(\epsilon))^{-1} p_z,
\end{equation} 

\noindent where $\eta = \eta(\epsilon) > 0$ is the unique solution of the equation

\begin{equation*}
\eta^2 \sum_{j=1}^m \frac{(\lambda_j - \epsilon^2)\|P_j p_z\|^2}{(\lambda_j + \eta)^2} = \epsilon^2 \|P_0 p_z\|^2.
\end{equation*}

\noindent In particular, $\eta(\epsilon) \sim \epsilon / \sigma$ as $\epsilon \to 0$.

\end{theorem}

The proof of this theorem is based on \cite{grho-annulus, grho-CEMP} with a few technical differences, so for completeness we will present the full argument here. Let us make a few remarks before presenting the proof:

\begin{enumerate}

\item In the formula \eqref{A_z using u_eps} for $A_z(\epsilon)$ we can take limits as $\epsilon \to 0$ and obtain (see Section~\ref{SECT Proof})

\begin{equation*}
A_z(0) = \frac{P_0 p_z(z)}{\|P_0 p_z\|} = \|P_0 p_z\|.
\end{equation*}

In particular, the maximizer function for $A_z(0)$ is $P_0 p_z / \|P_0 p_z\|$ (note that the assumption $p_z \notin \PK(\PH)$ implies $P_0p_z \neq 0$).

\item  Let $G \in \mathbb{C}^{n \times n}$ be the Gram matrix with $G_{jk} = (p_{z_j}, p_{z_k}) = p_{z_j}(z_k)$, then $\CK f = \lambda f$ implies that $G^T \bm{f} = \lambda \bm{f}$. So $\{\lambda_j\}_1^m$ are also eigenvalues of $G$, and after finding the corresponsing eigenvectors $\{\bm{f}_j\}$, from \eqref{K} we see that the eigenfunctions of $\PK$ are linear combinations of the functions $\{p_{z_j}\}$ with coefficients given by the eigenvectors.

For example, in the trivial case when there is only one data point $z_1$ we have $\lambda_1 = \|p_{z_1}\|^2$, \ $P_0 = Id - P_1$ and 

$$P_1f = \frac{f(z_1)}{\|p_{z_1}\|^2} p_{z_1}.$$

\item If $p_z \in \PK (\PH)$, then $\PK u = p_z$ for some $u \in \PH$. This relation means that the value of any function at $z$ is determined by its values at $\{z_j\}$, indeed $f(z) = \sum_j c_j \overline{f(z_j)}$ for any $f \in \PH$ where we set $c_j = u(z_j)$. So in this case we have complete stability: order $\epsilon$ smallness on $\{z_j\}$ implies order $\epsilon$ smallness at $z$. Indeed, applying the Cauchy-Schwartz inequality to $f(z)$ and using that $[f]\leq \epsilon$ we obtain

$$A_z(\epsilon) \leq \epsilon |\bm{c}|, \qquad \qquad \bm{c}=(c_1,...,c_n). $$

Such situations arise if we consider spaces that contain "boundary conditions". For example, $S = \{f \in \PH: f(z) = f(z_1)\}$ is a closed subspace of $\PH$ and hence is a RKHS in its own right, to which the above discussion applies. 

\end{enumerate}


\section{Proof of Theorem~\ref{THM discrete opt rec}} \label{SECT Proof}
\setcounter{equation}{0}

First consider the case $p_z \in \ker \PK$. Let us show that $A_z(\epsilon) = A_z(0)$ for any $\epsilon$, and so \eqref{main identity} is satisfied with $\sigma = 0$. Indeed, using the Cauchy-Schwartz inequality

$$|(f, p_z)| \leq \|f\| \|p_z\| \leq \|p_z\|,$$

\noindent giving the trivial bound $A_z(\epsilon) \leq \|p_z\|$. But this yields an optimal bound in this case, since the function $f_* = p_z / \|p_z\|$ satisfies both of the constraints and $|f_*(z)| = \|p_z\|$, attains the upper bound.

So let us concentrate on the non-trivial case $p_z \notin \ker \PK$ (this assumption is used later in the proof, namely in \eqref{Phi(infty)}).

The two constraints of \eqref{A abstract discrete} are invariant under multiplying $f$ with a constant phase factor: if $f$ satisfies the constraints, then so does $\lambda f$ for any $\lambda \in \mathbb{C}$ with $|\lambda| = 1$. So instead of maximizing $|f(z)|$ we can equivalently maximize $\re f(z)$. Using the reproducing kernel property and \eqref{K quadratic form} we rewrite \eqref{A abstract discrete} as a convex maximization problem with a linear target functional and two quadratic constraints: 

\begin{equation} \label{max problem}
\begin{cases}
\re (f, p_z) \to \max
\\
(f, f) \leq 1
\\
(\CK f, f) \leq \epsilon^2
\end{cases}
\end{equation}

\noindent Introduce Lagrange multiplies: nonnegative numbers $\mu$ and $\nu$ such that $\mu + \nu \neq 0$. Multiply the first constraint of \eqref{max problem} by $\mu$, the second one by $\nu$ and add the two inequalities to obtain

\begin{equation} \label{Lagrange}
((\mu + \nu \PK)f, f) \le \mu + \nu \epsilon^{2}.
\end{equation}

\noindent Now, if $M$ is a uniformly positive definite self-adjoint operator on $\PH$, expanding $(M(M^{-1}g - f),(M^{-1}g - f))\geq 0$, we obtain that for any $f, g \in \PH$

\[
2 \re (f, g) - (M^{-1}g, g) \le (M f, f).
\]

\noindent The uniform positivity of $M$ ensures that $M^{-1}$ is defined on all of $\PH$. This is an example of convex duality (cf. \cite{ekte76}) applied to the convex function $f \mapsto (M f, f)/2$. Then we also have for $\mu>0$, taking $M=\mu+\nu \CK$ and $g=p_z$ in the above inequality we get

\begin{equation}\label{cxdual}
2 \Re(f, p_z) - \left((\mu+\nu \CK)^{-1} p_z, p_z \right) \le
\left( (\mu+\nu \CK) f, f \right) \le
\mu + \nu \epsilon^{2} ,
\end{equation}

\noindent so that

\begin{equation} \label{maxub}
\Re(f, p_z) \le \frac{1}{2} \left((\mu+\nu \CK)^{-1} p_z, p_z \right)
+\frac{1}{2} \left( \mu+\nu \epsilon^{2}  \right),
\end{equation}

\noindent which is valid for every $f$, satisfying the constraints of \eqref{max problem} and all $\mu > 0$, $\nu\ge
0$. In order for the bound to be optimal we must have equality in
(\ref{cxdual}), which holds iff $p_z = (\mu + \nu \CK) f$ giving the formula for optimal function $f$:

\begin{equation} \label{maxxi}
f = (\mu + \nu \CK)^{-1} p_z.
\end{equation}

\noindent The goal is to choose the Lagrange multipliers $\mu$ and $\nu$ so that the constraints in \eqref{max problem} are satisfied by $f$, given by (\ref{maxxi}).

\begin{enumerate}

\item if $\nu = 0$, then $f = \frac{p_z}{\|p_z\|}$ does not depend on the small parameter $\epsilon$, which leads to a contradiction, because $p_z \notin \ker \CK$ implies that $(\CK f, f)>0$ and hence the second constraint $(\CK f, f) \leq \epsilon^2$ is violated if $\epsilon$ is small enough.

\item if $\mu= 0$, then $\CK f = \frac{1}{\nu} p_z$, contradicting to our assumption $p_z \notin \PK(\PH)$.

\end{enumerate}

Thus we are looking for $\mu>0,\  \nu> 0$, so that equalities hold  at both of the constraints for the function \eqref{maxxi} (these are the complementary slackness relations in Karush-Kuhn-Tucker conditions.), i.e.

\begin{equation}\label{cxdeq}
\left\| (\mu + \nu \CK)^{-1} p_z \right\| = 1, \qquad \qquad
\left[ (\mu + \nu \CK)^{-1} p_z \right] = \epsilon.
\end{equation}

\noindent Let $\eta = \frac{\mu}{\nu}$, solving the first equation of \eqref{cxdeq} for $\nu$ we find $\nu = \|(\CK+\eta)^{-1} p_z\|$. Then the square of the second equation of \eqref{cxdeq} reads

\begin{equation}\label{Phi(eta) def and equation}
\Phi(\eta):= \frac{ \left[ (\CK + \eta)^{-1} p_z  \right]^2 }
{ \|(\CK+\eta)^{-1} p_z\|^2 } = \epsilon^{2}.
\end{equation}

Let us now use the spectral decomposition of $\CK$. Recall that $P_j P_k = 0$ if $j\neq k$ and

\begin{equation*}
Id = \sum_{j=0}^m P_j,
\qquad \text{and} \qquad 
\CK = \sum_{j=1}^m \lambda_j P_j,
\end{equation*} 

\noindent further we also have

\begin{equation*}
\CK + \eta = \eta P_0 + \sum_{j=1}^m (\lambda_j + \eta) P_j,
\qquad \qquad
(\CK + \eta)^{-1} = \eta^{-1} P_0 + \sum_{j=1}^m \frac{P_j}{\lambda_j + \eta}.
\end{equation*}

\noindent Then writing the numerator of $\Phi$ as the quadratic form of $\PK$ using \eqref{K quadratic form} we get

\begin{equation*}
\Phi(\eta) = \frac{\displaystyle \sum_{j=1}^m \frac{\lambda_j \|P_j p_z\|^2}{(\lambda_j + \eta)^2} }{\displaystyle \eta^{-2} \|P_0 p_z\|^2 + \sum_{j=1}^m \frac{\|P_j p_z\|^2}{(\lambda_j + \eta)^2}}.
\end{equation*}

Next our goal is to show that the equation \eqref{Phi(eta) def and equation} has a unique solution $\eta = \eta(\Ge) > 0$. Clearly, $\Phi(0^+) = 0$ and 

\begin{equation} \label{Phi(infty)}
\Phi(+\infty) = \frac{[p_z]^2}{\|p_z\|^2} = \frac{(\CK p_{z},p_{z})}{\|p_{z}\|^{2}} > 0,
\end{equation}

\noindent because $p_z \notin \ker \PK$. Therefore, showing that $\Phi$ is strictly increasing will imply that \eqref{Phi(eta) def and equation} has a unique solution for $\epsilon$ small enough, namely for any $\epsilon^2 < \Phi(+\infty)$. So let us prove that $\Phi'(\eta)>0$. Set $a_j = \|P_jp_z\|^2$ for $j=0,...,m$, then the numerator of $\Phi'(\eta)$ (up to a factor of 2) can be simplified to

\begin{equation} \label{numerator}
\frac{a_0}{\eta^3} \sum_{j} \frac{\lambda_j^2 a_j}{(\lambda_j+\eta)^2} + \sum_{j,k} \lambda_j b_{jk} (\lambda_j - \lambda_k),
\end{equation}

\noindent where $b_{jk} = a_j a_k / (\lambda_j+\eta)^3 (\lambda_k+\eta)^3$. In particular $b_{jk} = b_{kj}$, so splitting the second sum of \eqref{numerator} into parts where $j>k, \ j<k$ and swapping the indices $j,k$ in the second part we find

$$\sum_{j,k} \lambda_j b_{jk} (\lambda_j - \lambda_k) = \sum_{j>k} \lambda_j b_{jk} (\lambda_j - \lambda_k) + \sum_{j>k} \lambda_k b_{kj} (\lambda_k - \lambda_j) = \sum_{j>k} b_{jk} (\lambda_j - \lambda_k)^2.$$  

\noindent Recalling that all $\lambda_j$ are distinct and positive we conclude that \eqref{numerator} is strictly positive (the value $0$ is excluded since otherwise $a_j=0$ for all $j=1,...,m$ implying that $p_z \in \ker\PK$) and hence $\Phi'(\eta)>0$.

Observe that, $\Phi(\eta) \sim \sigma^2 \eta^2$ as $\eta \to 0$, where $\sigma$ is given by \eqref{sigma}. Hence for the solution of the equation \eqref{Phi(eta) def and equation} we have $\eta(\epsilon) \sim \epsilon / \sigma$ as $\epsilon \to 0$. 

Setting $u=(\CK + \eta(\epsilon))^{-1}p_z$, from \eqref{maxub} we obtain

\begin{equation*}
\re (f,p_z)\le \frac{(u,p_z)}{2\|u\|} + \frac{\|u\|}{2} (\epsilon^2 + \eta(\epsilon)).
\end{equation*}

\noindent Definitions of $u$ and $\eta(\epsilon)$ imply that $(\CK u, u) / \|u\|^2 = \epsilon^2$, on the other hand

\begin{equation*}
u(z) = (u, p_z) = (u, \CK u + \eta(\epsilon) u) = (\CK u, u) + \eta(\epsilon) \|u\|^2 = (\epsilon^2 + \eta(\epsilon)) \|u\|^2,
\end{equation*}

\noindent which implies the optimal bound

\[
|f(z)|=\Re(f,p_z) \le \frac{u(z)}{2\|u\|} + \frac{u(z)}{2\|u\|} = \frac{u(z)}{\|u\|}.
\]

\noindent Thus

\begin{equation*}
A_z(\epsilon) = \frac{u(z)}{\|u\|} = [\epsilon^2 + \eta(\epsilon)] \|u\|.
\end{equation*}

\noindent It remain to analyze the asymptotic behavior of $A_z(\epsilon)$. To that end, using the spectral decomposition of $\PK$, note that

\begin{equation*}
\|u\|^2 = \eta(\epsilon)^{-2} a_0 + \sum_{j=1}^m \frac{a_j}{(\lambda_j + \eta(\epsilon))^2},
\end{equation*}

\noindent but then

\begin{equation*}
A_z^2(\epsilon) = \left(1 + \frac{\epsilon^2}{\eta(\epsilon)} \right)^2 a_0 + [\epsilon^2 + \eta(\epsilon)]^2 \sum_{j=1}^m \frac{a_j}{(\lambda_j + \eta(\epsilon))^2}.
\end{equation*}

\noindent Letting $\epsilon \to 0$ in the last formula gives $A_z^2(0) = a_0$. Finally, using the asymptotics of $\eta(\epsilon)$ we conclude that

\begin{equation*}
A_z^2(\epsilon) = (1+2\sigma \epsilon) A_z^2(0) + O(\epsilon^2), 
\end{equation*}

\noindent which, upon using the expansion $\sqrt{1+x} = 1+ \frac{x}{2} + O(x^2)$ for small $x$, implies the relation \eqref{main identity} and concludes the proof.

\bibliographystyle{abbrv}
\bibliography{refs}

\def\cprime{$'$} \ifx \cedla \undefined \let \cedla = \c \fi\ifx \cyr
  \undefined \let \cyr = \relax \fi\ifx \cprime \undefined \def \cprime
  {$\mathsurround=0pt '$}\fi\ifx \prime \undefined \def \prime {'}
  \fi\def\Ya{Ya}
\begin{thebibliography}{10}

\bibitem{berg78}
D.~J. Bergman.
\newblock The dielectric constant of a composite material --- {A} problem in
  classical physics.
\newblock {\em Phys. Rep.}, 43:377--407, 1978.

\bibitem{brune31}
O.~Brune.
\newblock Synthesis of a finite two-terminal network whose driving-point
  impedance is a prescribed function of frequency.
\newblock {\em Journal of Mathematics and Physics}, 10(1-4):191--236, 1931.

\bibitem{capr74}
I.~Caprini.
\newblock On the best representation of scattering data by analytic functions
  in ${L}\sb{2}$-norm with positivity constraints.
\newblock {\em Nuovo Cimento A (11)}, 21:236--248, 1974.

\bibitem{capr79}
I.~Caprini.
\newblock Integral equations for the analytic extrapolation of scattering
  amplitudes with positivity constraints.
\newblock {\em Nuovo Cimento A (11)}, 49(3):307--325, 1979.

\bibitem{ciulli69}
S.~Ciulli.
\newblock A stable and convergent extrapolation procedure for the scattering
  amplitude.---i.
\newblock {\em Il Nuovo Cimento A (1965-1970)}, 61(4):787--816, Jun 1969.

\bibitem{davis52}
P.~Davis.
\newblock An application of doubly orthogonal functions to a problem of
  approximation in two regions.
\newblock {\em Transactions of the American Mathematical Society},
  72(1):104--137, 1952.

\bibitem{deto18}
L.~Demanet and A.~Townsend.
\newblock Stable extrapolation of analytic functions.
\newblock {\em Foundations of Computational Mathematics}, 19(2):297--331, 2018.

\bibitem{ekte76}
I.~Ekeland and R.~Temam.
\newblock {\em Convex analysis and variational problems}.
\newblock North-Holland Publishing Co., Amsterdam, 1976.
\newblock Translated from the French, Studies in Mathematics and its
  Applications, Vol. 1.

\bibitem{epst08}
C.~L. Epstein.
\newblock {\em Introduction to the mathematics of medical imaging}, volume 102.
\newblock Siam, 2008.

\bibitem{feyn64}
R.~P. Feynman, R.~B. Leighton, and M.~Sands.
\newblock {\em The {F}eynman lectures on physics. {V}ol. 2: {M}ainly
  electromagnetism and matter}.
\newblock Addison-Wesley Publishing Co., Inc., Reading, Mass.-London, 1964.

\bibitem{fisher83}
S.~Fisher.
\newblock {\em Function Theory on Planar Domains: A Second Course in Complex
  Analysis}.
\newblock A Wieley-Interscience publication. Wiley, 1983.

\bibitem{fran90}
J.~Franklin.
\newblock Analytic continuation by the fast fourier transform.
\newblock {\em SIAM journal on scientific and statistical computing},
  11(1):112--122, 1990.

\bibitem{fdfd09}
C.-L. Fu, Z.-L. Deng, X.-L. Feng, and F.-F. Dou.
\newblock A modified tikhonov regularization for stable analytic continuation.
\newblock {\em SIAM Journal on Numerical Analysis}, 47(4):2982--3000, 2009.

\bibitem{grs01}
B.~Girod, R.~Rabenstein, and A.~Stenger.
\newblock {\em Signals and systems}.
\newblock Wiley,, 2001.

\bibitem{grho-CEMP}
Y.~Grabovsky and N.~Hovsepyan.
\newblock On feasibility of extrapolation of the complex electromagnetic
  permittivity function using {K}ramer-{K}ronig relations.
\newblock submitted.

\bibitem{grho-annulus}
Y.~Grabovsky and N.~Hovsepyan.
\newblock Explicit power laws in analytic continuation problems via reproducing
  kernel hilbert spaces.
\newblock {\em Inverse Problems}, 36(3):035001, feb 2020.

\bibitem{grho-gen}
Y.~Grabovsky and N.~Hovsepyan.
\newblock Optimal error estimates for analytic continuation in the upper
  half-plane.
\newblock {\em Communications on Pure and Applied Mathematics}, 74(1):140--171,
  2021.

\bibitem{lali60:8}
L.~D. Landau and E.~M. Lifshitz.
\newblock {\em Electrodynamics of continuous media}, volume~8.
\newblock Pergamon, New York, 1960.
\newblock Translated from the Russian by J. B. Sykes and J. S. Bell.

\bibitem{lipt01a}
R.~Lipton.
\newblock An isoperimetric inequality for gradients of solutions of elliptic
  equations in divergence form with applicatuion to the design of two-phase
  heat conductors.
\newblock to appear in SIAM J. Math. Anal.

\bibitem{lipt01}
R.~Lipton.
\newblock Optimal inequalities for gradients of solutions of elliptic equations
  occuring in two-phase heat conductors.
\newblock preprint.

\bibitem{mafed01}
L.~Maergoiz and A.~Fedotov.
\newblock Optimal error of analytic continuation from a finite set with
  inaccurate data in hilbert spaces of holomorphic functions.
\newblock {\em Siberian Mathematical Journal}, 42:926--935, 2001.

\bibitem{maos75}
A.~Marchuk and K.~Osipenko.
\newblock Best approximation of functions specified with an error at a finite
  number of points.
\newblock {\em Mathematical Notes of the Academy of Sciences of the USSR},
  17:207–--212, 1975.

\bibitem{micriv77}
C.~A. Micchelli and T.~J. Rivlin.
\newblock {\em A Survey of Optimal Recovery}, pages 1--54.
\newblock Springer US, Boston, MA, 1977.

\bibitem{micriv85}
C.~A. Micchelli and T.~J. Rivlin.
\newblock Lectures on optimal recovery.
\newblock In P.~R. Turner, editor, {\em Numerical Analysis Lancaster 1984},
  pages 21--93, Berlin, Heidelberg, 1985. Springer Berlin Heidelberg.

\bibitem{mill70}
K.~Miller.
\newblock Least squares methods for ill-posed problems with a prescribed bound.
\newblock {\em SIAM Journal on Mathematical Analysis}, 1(1):52--74, 1970.

\bibitem{mi80}
G.~W. Milton.
\newblock Bounds on complex dielectric constant of a composite material.
\newblock {\em Appl. Phys. Lett.}, 37(3):300--302, 1980.

\bibitem{nik54}
H.~Nikaido.
\newblock On von neumann's min—max theorems.
\newblock {\em Pacific Journal of Mathematics}, 4:65–70, 1954.

\bibitem{para16}
V.~I. Paulsen and M.~Raghupathi.
\newblock {\em An introduction to the theory of reproducing kernel Hilbert
  spaces}, volume 152.
\newblock Cambridge University Press, 2016.

\bibitem{payne75}
L.~E. Payne.
\newblock {\em Improperly posed problems in partial differential equations},
  volume~22.
\newblock Siam, 1975.

\bibitem{trefe19}
L.~N. Trefethen.
\newblock Quantifying the ill-conditioning of analytic continuation.
\newblock {\em arXiv preprint arXiv:1908.11097}, 2019.

\bibitem{vese99}
S.~Vessella.
\newblock A continuous dependence result in the analytic continuation problem.
\newblock {\em Forum Mathematicum}, 11(6):695--703, 1999.

\end{thebibliography}
\end{document}